\documentclass{elsarticle}

\usepackage{amssymb}

\newtheorem{theorem}{Theorem}
\newtheorem{statement}{Statement}

\newtheorem{corollary}[theorem]{Corollary}

\newtheorem{example}[theorem]{Example}

\newtheorem{remark}[theorem]{Remark}

\newfont{\eurorm}{eurm10 scaled 1100}
\newfont{\eurorms}{eurm10 scaled 800}

\journal{arXiv.org}

\begin{document}

\begin{frontmatter}

\title{\bf \normalsize
 SHARPENING  AND GENERALIZATIONS OF SHAFER-FINK \\[0.25 ex] AND WILKER TYPE INEQUALITIES: A NEW APPROACH}

\author{Marija Ra\v sajski}
\ead{marija.rasajski@etf.rs}
\author{Tatjana Lutovac}
\ead{tatjana.lutovac@etf.rs}
\author{Branko Male\v sevi\' c${}^{\ast}$}
\ead{branko.malesevic@etf.rs}

\cortext[cor]{Corresponding author.                                                               \\
Research of the first and second
and  third author was supported in part by the Serbian Ministry of
Education, Science and Technological Development, under Projects
ON 174033, TR 32023 and ON 174032 \& III 4400, respectively.}

\address{\rm School of Electrical Engineering, University of  Belgrade,                           \\[-0.25 ex]
Bulevar kralja Aleksandra 73, 11000 Belgrade, Serbia                                              \\[-3.25 ex]
}

\begin{abstract}
In this paper we propose and prove some generalizations and sharpenings of certain inequalities
of {\sc Wilker}'s and {\sc Shafer}-{\sc Fink}'s type. Application of the {\sc Wu}-{\sc Debnath}
theorem enabled us to prove some double sided inequalities.
\end{abstract}

\begin{keyword}
Sharpening; Generalization; Inequalities of {\sc Wilker}'s and {\sc Shafer}-{\sc Fink}'s type

\MSC 33B10; 26D05
\end{keyword}

\end{frontmatter}


\section{Introduction}


The main topic of this paper is the refinement and generalization of some  inequalities of {\sc Wilker}'s and {\sc Shafer}-{\sc Fink}'s type.

\medskip
{\sc Wilker}'s inequality is an inequality of the following form:
\begin{equation}
\left(\frac{\sin x}{x}\right)^{\!2}\!+\frac{\tan x}{x} > 2,
\end{equation}
and holds for  $x \!\in\! \left(0,\mbox{\small $\displaystyle\frac{\pi}{2}$}\right)$, \cite{J_B_Wilker_1989}.

That inequality had a great impact on numerous papers that address the theory of analytical inequalities \cite{G_V_Milovanovic_2014}.

Concerning {\sc Wilker}'s inequality, in this paper we propose and prove some extensions of Theorem 2.1 from \cite{Nenezic2016}, see also  \cite{C_Mortici_2014}.

\medskip
{\sc Shafer}-{\sc Fink}'s inequality is the following double-sided inequality:
\begin{equation}
\frac{3x}{2+\sqrt{1-x^2}} < \mbox{\rm arcsin}\,x < \frac{\pi x}{2+\sqrt{1-x^2}},
\end{equation}
and holds for $x \!\in\! (0,1)$, \cite{D_S_Mitrinovic_1970}, \cite{A_M_Fink_1995}.

\break

The above-mentioned inequality also had a great impact on many papers in the theory of analytical inequalities \cite{G_V_Milovanovic_2014}.
It is important to say that inequalities of this type have applications in various fields of engineering \cite{deAbreu_2009}, \cite{deAbreu_2012}, \cite{Alirezaei_Mathar_2014}; see also \cite{CloudDrachmanLebedev_2014}, \cite{Malesevic2017a}, \cite{Malesevic2018}.

In this paper, concerning {\sc Shafer}-{\sc Fink}'s inequality, we propose and prove some extensions of Theorems 1 and 2 from \cite{MRL_2017}.

\medskip
We now state the {\sc Wu}-{\sc Debnath} theorem (Theorem 2 in \cite{Wu_Debnath_2009}), used in our proofs.

\medskip
{\bf Theorem WD.}
\label{Debnath_Wu_T}
{\em Suppose that $f(x)$ is a real function on $(a,b)$, and that $n$ is a positive integer such that $f ^{(k)}(a+), f^{(k)}(b-)$,
$\left(k \!\in\! \{0,1,2, \ldots ,n\}\right)$ exist.

\medskip
\noindent
{\boldmath $(i)$} Supposing that $(-1)^{(n)} f^{(n)}(x)$~is~in\-cre\-asing on $(a,b)$, then
for all $x \in (a,b)$ the following inequality holds$\,:$
\begin{equation}
\label{Debnath_Wu_first}
\begin{array}{c}
\displaystyle\sum_{k=0}^{n-1}{\mbox{\small $\displaystyle\frac{f^{(k)}(b\mbox{\footnotesize $-$})}{k!}$}(x\!-\!b)^k}
+
\frac{1}{(a-b)^n}
{\bigg (}\!
f(a\mbox{\footnotesize $+$})
-
\displaystyle\sum_{k=0}^{n-1}{\mbox{\small $\displaystyle\frac{(a\!-\!b)^{k}f^{(k)}(b\mbox{\footnotesize $-$})}{k!}$}
\!{\bigg )} (x\!-\!b)^{n}}                                                                                    \\[2.5 ex]
<
f(x)
<
\displaystyle\sum_{k=0}^{n}{\frac{f^{(k)}(b\mbox{\footnotesize $-$})}{k!}(x\!-\!b)^{k}}.
\end{array}
\end{equation}
Furthermore, if $(-1)^{n} f^{(n)}(x)$ is decreasing on $(a,b)$, then the reversed inequality of {\rm (\ref{Debnath_Wu_first})} holds.

\medskip
\noindent
{\boldmath $(ii)$} Supposing that $f^{(n)}(x)$ is increasing on $(a,b)$, then for all $x \!\in\! (a,b)$
the following inequality also holds$\,:$
\begin{equation}
\label{Debnath_Wu_second}
\begin{array}{c}
\displaystyle\sum_{k=0}^{n-1}{\mbox{\small $\displaystyle\frac{f^{(k)}(a\mbox{\footnotesize $+$})}{k!}$}(x\!-\!a)^k}
+
\frac{1}{(b\!-\!a)^n}
{\bigg (}\!
f(b-)
-
\displaystyle\sum_{k=0}^{n-1}{\mbox{\small $\displaystyle\frac{(b-a)^{k}f^{(k)}(a\mbox{\footnotesize $+$})}{k!}$}
\!{\bigg )} (x\!-\!a)^{n}}                                                                                    \\[2.5 ex]
>
f(x)
>
\displaystyle\sum_{k=0}^{n}{\frac{f^{(k)}(a\mbox{\footnotesize $+$})}{k!}(x-a)^{k}}.
\end{array}
\end{equation}
Furthermore, if $f^{(n)}(x)$ is decreasing on $(a,b)$, then the reversed inequality~of~\mbox{\rm (\ref{Debnath_Wu_second})} holds.}

\medskip
Let us mention that an interesting application of Theorem WD to the inequalities that involve hyperbolic functions was considered in \cite{Milica_Makragic_2017}.

\medskip
Here, we prove a statement that represents a natural extension of the above theorem.

\begin{theorem}
\label{Natural_Extension_Theorem}

For the function $f:(a,b) \longrightarrow R$ let there exist the power series expansion$\,:$
\begin{equation}
f(x)
=
\displaystyle\sum_{k=0}^{\infty}{c_{k}(x-a)^k},
\end{equation}
for every $x \!\in\! (a,b)$, where $\{c_{k}\}_{k \in N_0}$ is the sequence of coefficients such that there is only
a finite number of negative coefficients, and their indices are all in the set $J \!=\! \{j_0,\ldots,j_\ell\}$.

\noindent
Then, for the function
\begin{equation}
F(x)
=
f(x)-\displaystyle\sum_{i=0}^{\ell}{c_{j_i}(x-a)^{j_i}}
=
\displaystyle\sum_{k \in N_0 \backslash\!\;\! J}{c_{k}(x-a)^k},
\end{equation}
and the sequence $\{C_{k}\}_{k \in N_0}$ of the non-negative coefficients defined by
\begin{equation}
C_{k}
=
\left\{
\begin{array}{ccc}
c_{k} \!&\!:\!&\! c_{k} > 0,               \\[0.75 ex]
0     \!&\!:\!&\! c_{k} \leq 0;
\end{array}
\right.
\end{equation}
holds that$\,:$
\begin{equation}
F(x)
=
\displaystyle\sum_{k=0}^{\infty}{C_{k}(x-a)^k},
\end{equation}
for every $x \!\in\! (a,b)$.

\smallskip
It is also $F^{(k)}(a+)= k!C_{k}$ $(k \!\in\! \{0,1,$ $2, \ldots ,n\})$ and the following inequalities hold$\,:$
\begin{equation}
\!
\begin{array}{c}
\displaystyle\sum_{k=0}^{n-1}{C_k(x-a)^k}
+
\frac{1}{(b-a)^n}
{\bigg (}
F(b-)
-
\displaystyle\sum_{k=0}^{n-1}{(b-a)^kC_k}
{\bigg )} (x-a)^{n}                                                                                          \\[1.5 ex]
>
F(x)
>
\displaystyle\sum_{k=0}^{n}{C_k(x-a)^{k}},
\end{array}
\end{equation}
i.e.$\,:$
\begin{equation}
\!\!
\begin{array}{c}
\displaystyle\sum_{k = 0}^{n - 1}{\!C_k}{(x \!-\! a)^{k}}
\!+\!
\displaystyle\sum_{i = 0}^{\ell}{\!c_{j_i}}{(x \!-\! a)^{j_i}}
\!+\!
\displaystyle\frac{{(x \!-\! a)}^n}{{(b \!-\! a)}^n}
\!\left(\!{f(b\mbox{\footnotesize $-$})
\!-\!
\displaystyle\sum_{k=0}^{n-1}{\!C_k}{{(b\!-\!a)}^k}
\!-\!
\displaystyle\sum_{i=0}^{\ell}{\!c_{j_i}}{(b\!-\!a)}^{j_i}}\!\right)                                            \\[1.5 ex]
> f(x) >
\displaystyle\sum_{k=0}^{n}{\!C_k}{(x\!-\!a)^k}
+
\displaystyle\sum_{i=0}^{\ell}{\!c_{j_i}}{(x - a)^{j_i}},
\end{array}
\end{equation}

\vspace*{-1.0 mm}

\noindent
for every $x \!\in\! (a,b)$

\end{theorem}
{\bf Proof.} This is a direct consequence of the previous theorem. The fact that all coefficients are positive implies that
all derivatives are positive, and, consequently, corresponding functions are increasing. \hfill $\;\Box$

\begin{corollary}
Let there hold the conditions from the previous theorem. If
\begin{equation}
n > \max\{j_0,\ldots,j_\ell\},
\end{equation}
then the following holds$\,:$
\begin{equation}
\!
\begin{array}{c}
\displaystyle\sum_{k=0}^{n-1}{\!c_k(x-a)^k}
+
\frac{1}{(b-a)^n}
{\bigg (}
f(b-)
-
\displaystyle\sum_{k=0}^{n-1}{\!c_k(b-a)^k}
{\bigg )}(x-a)^{n}                                                                                            \\[1.5 ex]
>
f(x)
>
\displaystyle\sum_{k=0}^{n}{\!c_k(x-a)^{k}},
\end{array}
\end{equation}
for every $x \!\in\! (a,b)$.
\end{corollary}

\break


\section{Main results}

\subsection{Wilker's type inequalities}

\medskip
The following statement was proved in \cite{Nenezic2016}:
\begin{theorem}
$\big ($Theorem $2.1$, {\rm \cite{Nenezic2016}}$\big )$
For every $x \!\in\! \left(0, \mbox{\small $\displaystyle\frac{\pi}{2}$} \right)$ the following inequalities are true$\,:$
\begin{equation}
\label{NMM_2016}
2+\left(\mbox{\small $\displaystyle\frac{8}{45}$}-a\left(x\right) \right) x^3 \tan x
<
\left(\mbox{\small $\displaystyle\frac{\sin x}{x}$} \right)^{\!2}+\mbox{\small $\displaystyle\frac{\tan x}{x}$}
<
2+\left(\mbox{\small $\displaystyle\frac{8}{45}$}-b_{1} \left(x \right)\right) x^3 \tan x,
\end{equation}
where $a(x)= \mbox{\small $\displaystyle\frac{8}{945}$} x^2$, $b_{1}(x)= \mbox{\small $\displaystyle\frac{8}{945}$} x^2
- \mbox{\small $\displaystyle\frac{\mbox{\boldmath $a$} }{14175}$}x^4$
and
$
\mbox{\boldmath $a$}
\!=\!
\mbox{\small $\displaystyle\frac{480\pi^6-40320\pi^4+3628800}{\pi^8}$}
$
$
=
17.15041 \ldots \, .
$
\end{theorem}
\quad\enskip
Above theorem is also an extension of Theorem 1 from \cite{C_Mortici_2014}.

\medskip
Let us notice that the inequality  (\ref{NMM_2016}) could be stated  as an equivalent inequality of the following form:
\begin{equation}
\frac{8}{45}x^3 - \frac{8}{945}x^5
<
\frac{1}{x} + \frac{\sin 2x}{2x^2} - 2 \cot x
<
\frac{8}{45}x^3 - \frac{8}{945}x^5 + \frac{\mbox{{\boldmath $a$}}}{14175}x^7,
\end{equation}
for $x \!\in\! \left(0, \mbox{\small $\displaystyle\frac{\pi}{2}$} \right)$.\\

In this paper, we sharpen the previous double-sided inequality using Theorem WD.

\begin{theorem}
For the function

\begin{equation}
f(x) =\frac{1}{x} + \frac{\sin 2x}{2x^2} - 2 \cot x - \frac{{8{x^3}}}{{45}} + \frac{{8{x^5}}}{{945}},
\end{equation}
where $x \!\in\! \left(0, \mbox{\small $\displaystyle\frac{\pi}{2}$} \right)$, the following sequence of inequalities holds$\,:$

\begin{equation}
\label{Th_2_inequalies}
\!\displaystyle\sum_{k=0}^{m}{c_{k}x^{2k+1}}
<
f(x)
<
\displaystyle\sum_{k=0}^{m-1}{c_{k}x^{2k+1}}
+
{\bigg (}
\!f\!\left({\frac{\pi}{2}}\right)
-
\!
\displaystyle\sum_{k=0}^{m-1}{\!c_{k}\left(\mbox{\small $\displaystyle\frac{\pi}{2}$}\right)^{\!2k+1}}
{\bigg )}\! \left(\mbox{\small $\displaystyle\frac{2x}{\pi}$}\right)^{\!2m+1}\!\!,
\end{equation}
for $x \!\in\! \left(0, \mbox{\small $\displaystyle\frac{\pi}{2}$} \right)$ and $m \in N$ and  $c_0=c_1=c_2=0$ and for $k \geq 3:$
\begin{equation}
\label{Def_c_k}
c_k = \frac{2^{2k+2}\left((4k\!+\!6) |\mbox{\boldmath $B$}_{2k+2}| + (-1)^{k+1} \right)}{(2k+3)!},
\end{equation}
where $\mbox{\boldmath $B$}_{i}$ are {\sc Bernoulli}'s numbers.
\end{theorem}
{\bf Proof.}
First, let us recall some well-known series expansions:

$$
\sin 2x
=
\sum_{k=1}^{\infty}{\frac{(-1)^{k}2^{2k-1}}{(2k-1)!}x^{2k-1}} \qquad {\big (}x \!\in\! R{\big )}
$$
and
$$
\cot x
=
\frac{1}{x}
\sum_{k=0}^{\infty}{\frac{|\mbox{\boldmath $B$}_{2k}|2^{2k}}{(2k)!}x^{2k}} \qquad {\big (}x \!\in\! (-\pi,0) \cup (0, \pi){\big )}.
$$
If we define  $f(0)\!=\!0$, then we have {\sc Taylor}'s expansion of the function $f(x)$ for $x\!=\!0:$
$$
f(x)
=
\displaystyle\sum_{k=0}^{\infty}{c_{k} x^{2k+1}},
$$
where $c_0=c_1=c_2=0$ and for $k \geq 3:$
$$
c_k
=
\frac{2^{2k+2}\left((4k+6) |\mbox{\boldmath $B$}_{2k+2}| + (-1)^{k+1} \right)}{(2k+3)!}.
$$
The obtained {\sc Taylor}'s expansion of the function  $f(x)$ converges
for $x \!\in\! \left(0, \mbox{\small $\displaystyle\frac{\pi}{2}$} \right)$.

\medskip
Based on  (\ref{Def_c_k}), it is evident that for the sequence
$$
c_3 \!=\!  \mbox{\small $\displaystyle\frac{16}{14175}$},
c_4 \!=\!  \mbox{\small $\displaystyle\frac{8}{467775}$},
c_5 \!=\!  \mbox{\small $\displaystyle\frac{3184}{638512875}$},
\ldots
$$
holds that $c_{k} \!>\! 0$ for $k \!\geq\! 3$. Then, the function
\begin{equation}
f(x)
=
\frac{{16{x^7}}}{{14175}} + \frac{{8{x^9}}}{{467775}} + \frac{{3184{x^{11}}}}{{638512875}} + \frac{{272{x^{13}}}}{{638512875}} + \frac{{7264{x^{15}}}}{{162820783125}}
+
\ldots\,,
\end{equation}
for $x \!\in\! \left(0, \mbox{\small $\displaystyle\frac{\pi}{2}$} \right)$, satisfies the condition:
$$
f^{(n)}(x) \!>\! 0
$$
for $x \!\in\! \left(0, \mbox{\small $\displaystyle\frac{\pi}{2}$} \right)$ and $n\in N$.
Then, for $n\in N$, the functions $f^{(n)}(x)$  are all increasing for $x \!\in\! \left(0, \mbox{\small $\displaystyle\frac{\pi}{2}$}\right)$, and there exist values $f ^{(k)}(0\mbox{\small $+$}), f^{(k)}(\frac{\pi}{2}\mbox{\small $-$})$
for every \mbox{$k \!\in\! \{0,1,\ldots ,n\}$}.

\bigskip
The right-hand side of the inequality  (\ref{Th_2_inequalies}) is obtained using Theorem WD. \hfill $\Box$\\
\begin{example}
Now, let us show several examples of approximations of the function  $f(x)$ obtained for $m=3,4,5,6$
and $x \!\in\! \left(0, \mbox{\small $\displaystyle\frac{\pi}{2}$} \right)\!:$
\begin{itemize}

\item For $m=3$ we get the double-sided inequality that was proved in {\rm \cite{Nenezic2016}}$\,:$
$$
\frac{{16}}{{14175}}{x^7}
<
f(x)
<
\left(\frac{2}{\pi}\right)^{\!\!7}\!\left({\frac{2}{\pi} - \frac{{\pi^3}}{45} + \frac{{\pi ^5}}{3780}}\right)x^7,
$$
and, in this way, a new  proof of the results from the paper {\rm \cite{Nenezic2016}} was obtained.

\smallskip
For $m>3$ the results that follow are higher accuracies.

\item For $m=4$ we have$\,:$
$$
\frac{16\,x^7}{14175}
+
\frac{8\,x^9}{467775}
<
f(x)
<
\frac{16\,x^7}{14175}
+
\left(\frac{2}{\pi} \right)^{\!\!9}
\!\left(\frac{2}{\pi} - \frac{\pi^3}{45} + \frac{\pi^5}{3780} - \frac{\pi^7}{113400}\right)\!x^9
$$
\item For $m=5$ we have$\,:$
$$
\begin{array}{c}
\displaystyle
\frac{16\,x^7}{14175} + \frac{8\,x^9}{467775} + \frac{3184\,x^{11}}{638512875} \, <                           \\[3.0 ex]
<
f(x)
<                                                                                                             \\[1.0 ex]
< \, \displaystyle
\frac{16 \, x^7}{14175} + \frac{8\,x^9}{467775}
+
\left(\frac{2}{\pi}\right)^{\!\!11}
\!\!\left(\frac{2}{\pi} - \frac{\pi^3}{45} + \frac{\pi^5}{3780} - \frac{\pi^7}{113400} - \frac{\pi^9}{29937600}\right)x^{11}
\end{array}
$$
\item For $m=6$ we have$\,:$
$$
\begin{array}{c}
\displaystyle
\frac{16\,x^7}{14175}
+
\frac{8\,x^9}{467775}
+
\frac{3184\,x^{11}}{638512875}
+
\frac{272\,x^{13}}{638512875} <                                                                               \\[2.5 ex]
< f(x) <                                                                                                      \\[1.0 ex]
< \, \displaystyle
\frac{16\,x^7}{14175}
+
\frac{8\,x^9}{467775}
+
\frac{3184\,x^{11}}{638512875} \; +                                                                           \\[1.0 ex]
+
\displaystyle
\left(\frac{2}{\pi}\right)^{\!\!11}
\!\!\left(
\frac{2}{\pi}
-
\frac{\pi^3}{45}
+
\frac{\pi^5}{3780}
-
\frac{\pi^7}{113400}
-
\frac{\pi^9}{29937600}
-
\frac{199\pi^{11}}{81729648000}
\right)x^{13}
\end{array}
$$
\end{itemize}
\end{example}
\begin{remark}
Let us note that Theorem WD enables us to estimate the error of approximation.
The difference between the right-hand side and the left-hand side of
the double-sided inequality in the previous theorem can be represented
by the following function$\,:$
$$
\!R_{m}(x)
\!=\!
\left(\!\!{f\!\left({\frac{\pi }{2}}\right)
-
\mathop \sum \limits_{k=3}^{m}\frac{2^{2k+2}\!\left({(4k+6)\left|{{\mbox{\boldmath $B$}_{2k+2}}}\right|
+
{(-1)}^{k+1}}\right)}{(2k+3)!}{\left({\frac{\pi}{2}}\right)}^{\!2k+1}}\!\right)\!\!
\left(\!\frac{2x}{\pi}\!\right)^{\!\!2m+1}\!\!,
$$
for $x \!\in\! \left(0, \mbox{\small $\displaystyle\frac{\pi}{2}$} \right)$.

\noindent\quad\enskip
The maximum values of the above-mentioned difference in the interval
$\left(0, \mbox{\small $\displaystyle\frac{\pi}{2}$} \right)$, for $m=3,4,5,6$, are shown in the table below$\,:$
\begin{equation}
\begin{array}{|c|c|c|c|c|c|c|} \hline
m        \!&\! 3          \!&\! 4           \!&\! 5           \!&\! 6            \\ \hline
R_{m}(x) \!&\! 0.00191501 \!&\! 0.000919303 \!&\! 0.000202959 \!&\! 0.0000519655 \\ \hline
\end{array}
\end{equation}
\end{remark}
\smallskip


\subsection{Shafer-Fink's type inequalities}


\medskip
Let us start from the following assertions proved by {\sc Bercu} in \cite{GB_2017}.

\begin{statement}
\label{GB_statement_1} $($Theorem $1$, {\rm \cite{GB_2017}}$)$ For
every real number $0 \leq x\leq 1$, the following two-sided
inequality holds$\,:$
\begin{equation}
\frac{x^5}{180} + \frac{x^7}{189} \leq \mbox{\rm arcsin}\,x - \frac{3x}{2 + \sqrt{1 - x^2}} \leq \frac{\pi - 3}{2}.
\end{equation}
\end{statement}

\begin{statement}
\label{GB_statement_2} $($Theorem $3$, {\rm \cite{GB_2017}}$)$ For
every $x \!\in\! [0,1]$ on the left-hand side and  every $x \!\in\!
[0, 0.871433]$ on the right-hand side, the following inequalities
hold true$\,:$
\begin{equation}
\left(
1 - \frac{\pi}{3} \right) x
+
\left( \frac{1}{6}
-
\frac{\pi}{18} \right)x^3 \leq \mbox{\rm arcsin} \, x
- \frac{\pi x}{
2 + \sqrt{1 - x^2}} \leq \left( 1 - \frac{\pi}{3} \right)x.
\end{equation}
\end{statement}
\begin{statement}
\label{GB_statement_3} $($Theorem $2$, {\rm \cite{GB_2017}}$)$ For
every $0 \leq x \leq 1$, one has$\,:$
\begin{equation}
\mbox{\rm arcsin}\,x - \frac{3x}{2 + \sqrt{1 - x^2}} \geq
\frac{a(x)}{2 + \sqrt{1 - x^2}},
\end{equation}
where $a(x) = (1/60)x^5 + (11/840)x^7$.
\end{statement}

\medskip
In \cite{MRL_2017} the authors proved the following theorem.

\begin{statement} $($Theorem $1$, {\rm \cite{MRL_2017}}$)$
For $x \!\in\! [0,1]$, $n \!\in\! N$ and  $k = 3 \vee k = \pi $ the following inequality holds$\,:$
\begin{equation}
\sum_{m=0}^{n}{D_{k}(m) x^{2m+1}}
\leq
\mbox{\rm arcsin}\,x
-
\displaystyle\frac{kx}{2+\sqrt{1-x^2}},
\end{equation}
where
\begin{equation}
\label{Dk_seq}
D_{k}(m)=
\displaystyle\frac{(2m)!}{(m!)^2 (2m\!+\!1)2^{2m}}
-
\!\left(\!
\frac{(-1)^mk}{3^{m+1}}
+\!
\displaystyle\sum_{i = 0}^{m - 1}{
\frac{k{{(-1)}^{m-1-i}}(2i)!}{3^{m-i}i!(i\!+\!1)!2^{2i + 1}}}
\!\right)\!>0,
\end{equation}
for $m \!\in\! N_0$, $m \!\geq\! 2$ ${\big (}D_{k}(0)\!=\!D_{k}(1)\!=\!0{\big )}$.
\end{statement}
\begin{remark}
For $n\!=\!3,\,k\!=\!3$ and $n\!=\!1,\,k\!=\!\pi$ we get the left-hand sides of the
inequalities stated in Theorems $1$ and $3$ from {\rm \cite{GB_2017}} by {\sc G. Bercu}.
\end{remark}

Now, let us consider the functions:
\begin{equation}
f_{k}(x) = \mbox{\rm arcsin}\,x - \displaystyle\frac{kx}{2+\sqrt{1-x^2}},
\end{equation}
for $x \!\in\! [0,1]$, and $k \!=\! 3 \vee k \!=\! \pi $. Then, using Theorem WD, we get:
\begin{theorem}

For $x \!\in\! [0,1]$ and the sequence $\{D_{k}(m)\}_{m \in N_0, m \geq 2}$ defined by {\em (\ref{Dk_seq})},
the following double-sided inequalities hold true$\,:$
\begin{equation}
\label{Ineq_D}
\!\!
\displaystyle\sum_{m=0}^{n}{\!D_k(m){x^{2m+1}}
<
f_k(x)
<
\displaystyle\sum_{m=0}^{n-1}
{\!D_k(m)x^{2m+1}}
\!+\!
\left(\!\,\!f_k(1)
\!-\!
\displaystyle\sum_{m=0}^{n-1}
D_k(m)\!\!\right)\!x^{2n+1}}.
\end{equation}
\end{theorem}
\begin{example}
We show a few examples of approximations of the function $f_{k}(x)$ for \mbox{$k=3$}, $n=3,4,5,6$
and $x \!\in\! \left(0, \mbox{\small $\displaystyle\frac{\pi}{2}$} \right)\!:$
\begin{itemize}
\item For $n=3$ we have$\,:$
$$
\frac{x^5}{180} + \frac{x^7}{189}
< f_3(x) <
\frac{x^5}{180} + \left( - \frac{271}{180} + \frac{\pi }{2} \right)\!x^7
$$
\item For $n=4$ we have$\,:$
$$
\frac{x^5}{180} + \frac{x^7}{189} + \frac{23x^9}{5184}
< f_3(x) <
\frac{x^5}{180} + \frac{x^7}{189} + \left( -\frac{5711}{3780} + \frac{\pi }{2} \right)\!x^9
$$
\item For $n=5$ we have$\,:$
$$
\frac{x^5}{180} + \frac{x^7}{189} + \frac{23x^9}{5184} + \frac{629x^{11}}{171072}
< f_3(x) <
\frac{x^5}{180} + \frac{x^7}{189} + \frac{23x^9}{5184} + \left( -\frac{274933}{181440} + \frac{\pi }{2} \right)\!x^{11}
$$
\item For $n=6$ we have$\,:$
$$
\begin{array}{c}
\displaystyle
\frac{x^5}{180} + \frac{x^7}{189} + \frac{23x^9}{5184} + \frac{629x^{11}}{171072} + \frac{14929x^{13}}{4852224} \, <  \\[2.5 ex]
< f_3(x) <                                                                                                            \\[1.0 ex]
< \, \displaystyle
\frac{x^5}{180} + \frac{x^7}{189} + \frac{23x^9}{5184} + \frac{629x^{11}}{171072}
+ \left( - \frac{2273701}{1496880}  + \frac{\pi }{2} \right){x^{13}}
\end{array}
$$

\end{itemize}
\end{example}
\begin{example}
Now, let us present several examples of approximations of the function $f_{k}(x)$ for \mbox{$k=\pi$}, $n=3,4,5,6$
and $x \!\in\! \left(0, \mbox{\small $\displaystyle\frac{\pi}{2}$} \right)\!:$
\begin{itemize}
\item  For $n=3$ we have$\,:$
$$
\begin{array}{c}
\left({1 \!-\! \frac{\pi }{3}} \right)\!x \!+\! \left( {\frac{1}{6} \!-\! \frac{\pi }{{18}}} \right)\!{x^3} \!+\! \left( {\frac{3}{{40}} \!-\! \frac{{5\pi }}{{216}}} \right)\!{x^5} \!+\! \left( {\frac{5}{{112}} \!-\! \frac{{17\pi }}{{1296}}} \right)\!{x^7} < f_{\pi}(x) < \\
\\
 < \left( {1 \!-\! \frac{\pi }{3}} \right)\!x \!+\! \left( {\frac{1}{6} \!-\! \frac{\pi }{{18}}} \right)\!{x^3} \!+\! \left( {\frac{3}{{40}} \!-\! \frac{{5\pi }}{{216}}} \right)\!{x^5} \!+\! \left( { \!-\! \frac{{149}}{{120}} \!+\! \frac{{89\pi }}{{216}}} \right)\!{x^7}
\end{array}
$$
\item For $n=4$ we have$\,:$
$$
\!\!\!\!\!\!\!
\!\!\!\!\!\!\!
\begin{array}{c}
\left({1 \!-\! \frac{\pi }{3}} \right)\!x \!+\! \left( {\frac{1}{6} \!-\! \frac{\pi }{{18}}} \right)\!{x^3} \!+\! \left( {\frac{3}{{40}} \!-\! \frac{{5\pi }}{{216}}} \right)\!{x^5} \!+\! \left( {\frac{5}{{112}} \!-\! \frac{{17\pi }}{{1296}}} \right)\!{x^7} \!+\! \left( {\frac{{35}}{{1152}} \!-\! \frac{{269\pi }}{{31104}}} \right)\!{x^9} < f_{\pi}(x) < \\
\\
 < \left( {1 - \frac{\pi }{3}} \right)\!x + \left( {\frac{1}{6} - \frac{\pi }{{18}}} \right)\!{x^3} + \left( {\frac{3}{{40}} - \frac{{5\pi }}{{216}}} \right)\!{x^5} + \left( {\frac{5}{{112}} - \frac{{17\pi }}{{1296}}} \right)\!{x^7} + \left( { - \frac{{2161}}{{1680}} + \frac{{551\pi }}{{1296}}} \right)\!{x^9}
\end{array}
$$
\item  For $n=5$ we have$\,:$
$$
\!\!\!\!\!\!\!\!\!\!\!\!\!\!\!\!\!\!\!\!
\!\!\!\!\!\!\!\!\!\!\!\!\!\!\!\!\!\!\!\!
\begin{array}{c}
\left( {1 \!-\! \frac{\pi }{3}} \right)\!x \!+\! \left( {\frac{1}{6} \!-\! \frac{\pi }{{18}}} \right)\!{x^3} \!+\! \left( {\frac{3}{{40}} \!-\! \frac{{5\pi }}{{216}}} \right)\!{x^5}
\!+\! \left( {\frac{5}{{112}} \!-\! \frac{{17\pi }}{{1296}}} \right)\!{x^7} \!+\! \left( {\frac{{35}}{{1152}} \!-\! \frac{{269\pi }}{{31104}}} \right)\!{x^9}
\!+\! \left( {\frac{{63}}{{2816}} \!-\! \frac{{1163\pi }}{{186624}}} \right)\!{x^{11}} < f_{\pi}(x) <                                                                                       \\[2.0 ex]
 < \left( {1 \!-\! \frac{\pi }{3}} \right)\!x \!+\! \left( {\frac{1}{6} \!-\! \frac{\pi }{{18}}} \right)\!{x^3} \!+\! \left( {\frac{3}{{40}} \!-\! \frac{{5\pi }}{{216}}} \right)\!{x^5}
 \!+\! \left( {\frac{5}{{112}} \!-\! \frac{{17\pi }}{{1296}}} \right)\!{x^7} \!+\! \left( {\frac{{35}}{{1152}} \!-\! \frac{{269\pi }}{{31104}}} \right)\!{x^9}
 \!+\! \left( { \!-\! \frac{{53089}}{{40320}} \!+\! \frac{{13493\pi}}{{31104}}} \right)\!{x^{11}}
\end{array}
$$
\item  For $n=6$ we have$\,:$
$$
\!\!\!\!\!\!\!\!\!\!\!\!\!\!\!\!\!\!\!\!\!\!\!\!\!\!\!\!\!\!\!
\!\!\!\!\!\!\!\!\!\!\!\!\!\!\!\!\!\!\!\!\!\!\!\!\!\!\!\!\!\!\!
\begin{array}{c}
\left( {1\mbox{\scriptsize $-$}\frac{\pi }{3}} \right)\!x\mbox{\footnotesize $+$}\left( {\frac{1}{6}\mbox{\scriptsize $-$}\frac{\pi }{{18}}} \right)\!{x^3}\mbox{\footnotesize $+$}\left( {\frac{3}{{40}}\mbox{\scriptsize $-$}\frac{{5\pi }}{{216}}} \right)\!{x^5}\mbox{\footnotesize $+$}\left( {\frac{5}{{112}}\mbox{\scriptsize $-$}\frac{{17\pi }}{{1296}}} \right)\!{x^7}\mbox{\footnotesize $+$}\left( {\frac{{35}}{{1152}}\mbox{\scriptsize $-$}\frac{{269\pi }}{{31104}}} \right)\!{x^9}\mbox{\footnotesize $+$}\left( {\frac{{63}}{{2816}}\mbox{\scriptsize $-$}\frac{{1163\pi }}{{186624}}} \right)\!{x^{11}}\mbox{\footnotesize $+$}\left( {\frac{{231}}{{13312}}\mbox{\scriptsize $-$}\frac{{10657\pi }}{{2239488}}} \right)\!{x^{13}}
< f_{\pi}(x) <                                                                                                          \\[2.0 ex]
< \left( {1\mbox{\scriptsize $-$}\frac{\pi }{3}} \right)\!x\mbox{\footnotesize $+$}\left( {\frac{1}{6}\mbox{\scriptsize $-$}\frac{\pi }{{18}}} \right)\!{x^3}\mbox{\footnotesize $+$}\left( {\frac{3}{{40}}\mbox{\scriptsize $-$}\frac{{5\pi }}{{216}}} \right)\!{x^5}\mbox{\footnotesize $+$}\left( {\frac{5}{{112}}\mbox{\scriptsize $-$}\frac{{17\pi }}{{1296}}} \right)\!{x^7}\mbox{\footnotesize $+$}\left( {\frac{{35}}{{1152}}\mbox{\scriptsize $-$}\frac{{269\pi }}{{31104}}} \right)\!{x^9}\mbox{\footnotesize $+$}\left( {\frac{{63}}{{2816}}\mbox{\scriptsize $-$}\frac{{1163\pi }}{{186624}}} \right)\!{x^{11}}\mbox{\footnotesize $+$}\left( {\mbox{\scriptsize $-$} \frac{{1187803}}{{887040}}\mbox{\footnotesize $+$}\frac{{82121\pi }}{{186624}}} \right)\!{x^{13}}
\end{array}
$$
\end{itemize}
\end{example}


\medskip

In \cite{MRL_2017} the authors proved the following assertion.

\begin{statement}
\label{improvement-Bercu-th-2}
$($Theorem $2$, {\rm \cite{MRL_2017}}$)$ If  $n \in N$ and  $n \geq 2$, then
\begin{equation}
\label{improvement3}
\mbox{\rm arcsin}\,x - \frac{3x}{2 +\sqrt{1-x^2}} \geq \frac{\mbox{\small
$\displaystyle\sum_{m=2}^{n}$}{E(m)x^{2m+1}}}{2 + \sqrt{1-x^2}},
\end{equation}
for every $x \!\in\! [0,1]$, where
\begin{equation}
\label{Em_seq}
E(m) = \frac{m \, (2m\!-\!1)!}{(2m\!+\!1)2^{2m-2} m!^2}
-
\frac{2m \, 2^{2m-2}(m\!-\!1)!^2}{(2m\!+\!1)!}>0,
\end{equation}
for $m \!\in\! N$, $m \!\geq\! 2$ ${\big (}E(1)\!=\!0{\big )}$.
\end{statement}
\begin{remark}
For $n\!=\!3$ we get the left-hand sides of the inequality stated in Theorem $2$ from {\rm \cite{GB_2017}} by {\sc G. Bercu}.
\end{remark}

Using Theorem WD we prove the following theorem.\\

\begin{theorem}

For $x \!\in\! [0,1]$ and the sequence $\{E(m)\}_{m \in N, m\geq 2}$, defined by {\rm (\ref{Em_seq})}, the following double-sided inequalities hold true$\,:$
\begin{equation}
\label{Ineq_E}
\mbox{\small $\displaystyle\frac{\displaystyle\sum_{m=2}^{n}{\!E(m)x^{2m+1}}}{2 +\sqrt{1-x^2}}$}
<
\mbox{\rm arcsin}\,x - \frac{3x}{2 +\sqrt{1-x^2}}
<
\mbox{\small $\displaystyle\frac{\displaystyle\sum_{m=2}^{n-1}{\!E(m)x^{2m+1}}
\!+\!
\left(\!\pi\!-\!
\displaystyle\sum_{m=0}^{n-1}{E(m)}\!\!\right)\!x^{2n+1}}{2 +\sqrt{1-x^2}}$},
\end{equation}
\end{theorem}
\begin{example}
Following are several examples of approximations of the function
\[
\mbox{\rm arcsin}\,x - \mbox{\small $\displaystyle\frac{3x}{2\!+\!\sqrt{1\!-\!x^2}}$}
\]
for $n=3,4,5,6$ and $x \!\in\! \left(0, \mbox{\small $\displaystyle\frac{\pi}{2}$} \right)\!:$

\begin{itemize}
\item For $n=3$ we have$\,:$
$$
\displaystyle
\frac{
\mbox{\small $\displaystyle\frac{1}{60}$}x^5
+
\mbox{\small $\displaystyle\frac{11}{840}$}x^7}{2 +\sqrt{1-x^2}}
<\mbox{\rm arcsin}\,x - \displaystyle\frac{3x}{2 +\sqrt{1-x^2}} <
\displaystyle
\frac{
\mbox{\small $\displaystyle\frac{1}{60}$}x^5
+
\mbox{\small $\displaystyle\left(\!\pi-\frac{181}{60}\!\right)$}x^{7}}{2 +\sqrt{1-x^2}}.
$$
\item For $n=4$ we have$\,:$
$$
\!\!\!\!\!\!\!
\displaystyle
\frac{
\mbox{\small $\displaystyle\frac{1}{60}$}x^5
\!+\!
\mbox{\small $\displaystyle\frac{11}{840}$}x^7
\!+\!
\mbox{\small $\displaystyle\frac{67}{6720}$}x^9}{2 +\sqrt{1-x^2}}
<\mbox{\rm arcsin}\,x - \displaystyle\frac{3x}{2 +\sqrt{1-x^2}} <
\displaystyle
\frac{
\mbox{\small $\displaystyle\frac{1}{60}$}x^5
\!+\!
\mbox{\small $\displaystyle\frac{11}{840}$}x^7
\!+\!
\mbox{\small $\displaystyle\left(\!\pi-\frac{509}{168}\!\right)$}x^{9}}{2 +\sqrt{1-x^2}}.
$$
\item For $n=5$ we have$\,:$
$$
\begin{array}{c}
\displaystyle
\frac{
\mbox{\small $\displaystyle\frac{1}{60}$}x^5
+
\mbox{\small $\displaystyle\frac{11}{840}$}x^7
+
\mbox{\small $\displaystyle\frac{67}{6720}$}x^9
+
\mbox{\small $\displaystyle\frac{3461}{443520}$}x^{11}}{2 +\sqrt{1-x^2}} \, <                                 \\[1.5 ex]
<\mbox{\rm arcsin}\,x - \displaystyle\frac{3x}{2 +\sqrt{1-x^2}}<                                              \\[1.5 ex]
\displaystyle
< \, \frac{
\mbox{\small $\displaystyle\frac{1}{60}$}x^5
+
\mbox{\small $\displaystyle\frac{11}{840}$}x^7
+
\mbox{\small $\displaystyle\frac{67}{6720}$}x^9
+
\mbox{\small $\displaystyle\left(\!\pi-\frac{6809}{2240}\!\right)$}x^{11}}{2 +\sqrt{1-x^2}}.
\end{array}
$$
\item For $n=6$ we have$\,:$
$$
\begin{array}{c}
\displaystyle
\frac{
\mbox{\small $\displaystyle\frac{1}{60}$}x^5
+
\mbox{\small $\displaystyle\frac{11}{840}$}x^7
+
\mbox{\small $\displaystyle\frac{67}{6720}$}x^9
+
\mbox{\small $\displaystyle\frac{3461}{443520}$}x^{11}
+
\mbox{\small $\displaystyle\frac{29011}{4612608}$}x^{13}}{2 +\sqrt{1-x^2}} \, <                               \\[1.5 ex]
< \mbox{\rm arcsin}\,x - \displaystyle\frac{3x}{2 +\sqrt{1-x^2}} <                                            \\[1.5 ex]
< \, \displaystyle
\frac{
\mbox{\small $\displaystyle\frac{1}{60}$}x^5
+
\mbox{\small $\displaystyle\frac{11}{840}$}x^7
+
\mbox{\small $\displaystyle\frac{67}{6720}$}x^9
+
\mbox{\small $\displaystyle\frac{3461}{443520}$}x^{11}
+
\mbox{\small $\displaystyle\left(\!\pi-\frac{1351643}{443520}\!\right)$}x^{13}}{2 +\sqrt{1-x^2}}.
\end{array}
$$
\end{itemize}
\end{example}

\section{Conclusion}

In this paper, we proposed and proved new inequalities, which represent refinements and generalizations of the inequalities
stated in \cite{MRL_2017} related to {\sc Shafer}-{\sc Fink}'s inequality for the inverse sine function, as well as the inequalities
stated in \cite{Nenezic2016} related to {\sc Wilker}'s inequality. Finally, let us note that proofs of inequalities
(\ref{Th_2_inequalies}) for any fixed $m \!\in\! N$, and inequalities (\ref{Ineq_D}) and (\ref{Ineq_E}) for any fixed
$n \!\in\! N$, can be obtained by methods and algorithms developed in \cite{Malesevic2016} and \cite{Lutovac2017}.



\smallskip
\noindent \textbf{Competing Interests.} The authors would like to
state that they do not have any competing interests in the subject
of this research.

\smallskip
\noindent \textbf{Author's Contributions.} All the authors
participated in every phase of the research conducted for this
paper.

\break


\break

\end{document}